\documentclass{llncs}
\usepackage{amssymb}
\usepackage{amsmath}
\usepackage{epsfig}
\usepackage{amsfonts}
\usepackage{latexsym}
\newcommand{\IR} {\mathbb{R}}
\newcommand{\IN} {\mathbb{N}}
\newcommand{\ZZ} {\mathbb{Z}}

\title{Discrete Tomography: Reconstruction under periodicity constraints}
\author{Alberto Del Lungo\inst{1} \and Andrea Frosini\inst{1}
\and Maurice Nivat\inst{2} \and Laurent Vuillon\inst{2}}
\institute{Dipartimento di Matematica, Universit\`a di Siena, Via
del Capitano 15, 53100, Siena, Italy
\email{[dellungo,frosini]@unisi.it} \and Laboratoire
d'Informatique, Algorithmique, Fondements et Applications (LIAFA)
Universit\'e Denis Diderot 2, place Jussieu 75251 Paris Cedex 05,
France \email{[Maurice.Nivat,Laurent.Vuillon]@liafa.jussieu.fr}}

\begin{document}
\maketitle

\begin{abstract}
This paper studies the problem of reconstructing binary matrices
that are only accessible through few evaluations of their discrete
X-rays. Such question is prominently motivated by the demand in
material science for developing a tool for the reconstruction of
crystalline structures from their images obtained by
high-resolution transmission electron microscopy. Various
approaches have been suggested for solving the general problem of
reconstructing binary matrices that are given by their discrete
X-rays in a number of directions, but more work have to be done to
handle the ill-posedness of the problem. We can tackle this
ill-posedness by limiting the set of possible solutions, by using
appropriate a priori information, to only those which are
reasonably typical of the class of matrices which contains the
unknown matrix that we wish to reconstruct. Mathematically, this
information is modelled in terms of a class of binary matrices to
which the solution must belong. Several papers study the problem
on classes of binary matrices on which some connectivity and
convexity constraints are imposed.

We study the reconstruction problem on some new classes consisting
of binary matrices with periodicity properties, and we propose a
polynomial-time algorithm for reconstructing these binary matrices
from their orthogonal discrete X-rays.

\vspace{0.2cm} \noindent {{\bf keywords:}\ combinatorial problem,
discrete tomography, binary matrix, polyomino, periodic
constraint, discrete X-rays.}
\end{abstract}

\section{Introduction}

The present paper studies the possibility of determining the
geometrical aspects of a discrete physical structure whose
interior is accessible only through a small number of measurements
of the atoms lying along a fixed set of directions. This is the
central theme of {\em Discrete Tomography} and the principal
motivation of this study is in the attempt to reconstruct
three-dimensional crystals from two-dimensional images taken by a
transmission electron microscope. The quantitative analysis of
these images can be used to determine the number of atoms in
atomic lines in certain directions~\cite{KSB,KSB2}. The question
is to deduce the local atomic structure of the crystal from the
atomic line count data. The goal is to use the reconstruction
technique for quality control in VLSI (Very Large Scale
Integration) technology. Before showing the results of this paper,
we give a brief survey of the relevant contributions in Discrete
Tomography.

Clearly, the best known and most important part of the general
area of tomography is {\em Computerized Tomography}, an invaluable
tool in medical diagnosis and many other areas including biology,
chemistry and material science. Computerized Tomography is the
process of obtaining the density distribution within a physical
structure from multiple X-rays. More formally, we attempts to
reconstruct a density function $f(x)$ for $x$ in $\IR^2$ or
$\IR^3$ from knowledge of its line integral $X_f(L)=\int_L f(x)
dx$ for lines $L$ through the space. This line integral is the
{\em X-ray} of $f(x)$ along $L$. The mapping $f\rightarrow X_f$ is
known as the {\em Radon transform}. The mathematics of
Computerized Tomography is quite well understood. Appropriate
quadratures~\cite{Lo,SK} of the Radon inversion formula are used,
with concepts from calculus and continuous mathematics playing the
main role.

Usually, the physical structure has a very big variety of density
values, and so a large number of  X-rays are necessary to ensure
the accurate reconstruction of their distribution. In some cases
the structure that we want to reconstruct has only a small number
of possible values. For example, a large number of objects
encountered in industrial computerized tomography (for the purpose
of non-destructive testing or reverse engineering) \cite{BKS} are
made of a single homogenous material. In many of these
applications there are strong technical reasons why only a few
X-rays of the structure can be physically determined. {\em
Discrete Tomography} is the area of Computerized Tomography in
which these special cases are studied. The name Discrete
Tomography is due to Larry Shepp, who organized the first meeting
devoted to the topics in 1994.

\begin{figure}[htb]
  \centerline{\epsfxsize=10.5cm \epsfbox{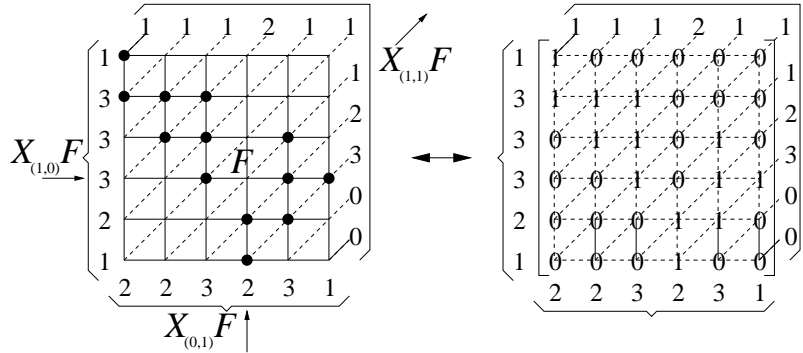} }
  \caption{\baselineskip=0ex \small
           A subset $F$ of $\ZZ^2$ with the corresponding binary matrix.
           $X_{(1,0)}F, X_{(0,1)}F$ and $X_{(1,1)}F$ are the discrete
           X-rays in the directions $(1,0),\;(0,1)$ and $(1,1)$.}
  \label{fig1}
\end{figure}

\noindent An example of such a case is the above-mentioned problem
of determining local atomic structure of a crystal from the atomic
line count data. In a simple but highly relevant model suggested
by Peter Schwander and Larry Shepp the possible atom locations in
a unit cell of a crystal are defined on the integer lattice
$\ZZ^3$, while the electron beams are modeled as lines parallel to
a given direction. The presence of an atom at a specific location
corresponds to a pixel value one at the location; the absence of
an atom corresponds to a pixel value zero. The number of atoms
along certain lines through the crystal (i.e., the sum of pixel
values along those lines) define the {\em discrete X-rays} of the
atomic structure with respect to the family of lines involved.
Since in practice, one degree of freedom for moving the imaging
device is used to control the position of the crystal, the view
directions for which data are provided lie all in the same plane.
This means that the 3D-problem leads itself to a 2D-slice-by-slice
reconstruction. Therefore, the crystal is represented by a binary
matrix and its discrete X-ray along a direction $u$ is an integral
vector giving the sum of its elements on each line parallel to $u$
(see Fig.~\ref{fig1}).

\noindent Measurements are usually only available along two, three
or four directions, which is much less than what is typical used
in Computerized Tomography (a few hundred). In fact, the electron
microscope makes measurements at the atomic level and uses high
energy (and so deeply penetrating) rays which can corrupt the
crystal itself. So, it can take only two, three or four images of
the crystal before the energy of the radiations destroys it, or at
least changes permanently its atomical configuration so that the
subsequent radiations will see something different from the
original one.

\noindent Now, the problem is to invert the {\em discrete Radon
transform}, i.e., to reconstruct the binary matrix from this small
set of discrete X-rays. More precisely, the basic question is to
determine, given a set of directions $u_1,\ldots, u_k$ and a set
of integral vectors $X_1, \ldots X_k$, whether there exists a
binary matrix $F$ whose discrete X-rays along $u_1,\ldots u_k$ are
$X_1,\ldots, X_k$. The general methods of Computerized Tomography
cannot be used effectively if the number of X-rays is so small,
and they seems unlikely to work in practice.

Discrete Tomography has its own mathematical theory mostly based
on discrete mathematics. It has some strong connection with
combinatorics and geometry. We wish to point out that the
mathematical techniques developed in Discrete Tomography have
applications in other fields such as: image processing~\cite{SH},
statistical data security~\cite{IJ}, biplane
angiography~\cite{PO}, graph theory~\cite{A} and so on. As a
survey of the state of the art of Discrete Tomography we can
suggest the book~\cite{HK}.

\noindent Interestingly, mathematicians have been concerned with
abstract formulations of these problems before the emergence of
the practical applications. Many problems of Discrete Tomography
were first discussed as combinatorial problems during the late
1950s and early 1960s. In 1957 Ryser~\cite{Ry} and Gale~\cite{Ga}
gave a necessary and sufficient condition for a pair of vectors
being the discrete X-rays of a binary matrix along  horizontal
and vertical directions. The discrete X-rays in horizontal and
vertical directions are equal to {\em row} and {\em column
sums} of the matrix. They gave an exact combinatorial
characterization of the row and column sums that correspond to a
binary matrix, and they derived a fast $O(nm)$ time algorithm for
reconstructing a matrix, where $n$ and $m$ denote its sizes. We
refer the reader to an excellent survey on the binary matrices
with given row and column sums by Brualdi~\cite{Br80}.

\noindent The space of solutions of the reconstruction problem,
however, is really huge and in general quite impossible to
control. A good idea may seem to start increasing the number of
X-rays one by one in order to decrease the number of
solutions. Unfortunately, the reconstruction problem becomes
intractable when the number of X-rays is greater than two, as
proved in \cite{GGP}. This means that (unless $P=NP$) exact
reconstructions require, in general, an exponential amount of
time. In polynomial time only approximate solutions can be
expected. In this context, an approximate solution is close to the
optimal one if its discrete X-rays in the set of prescribed
directions are close to those of the original set. Various
approaches have been suggested for solving the
problem~\cite{FSSV,TH,VL}. Recently, an interesting
method~\cite{GVW} for finding an approximate solutions has been
proposed. Even though the reconstruction problem is intractable,
some simple algorithms proposed in \cite{GVW} have good worst-case
bounds and they perform even better in computational practice.

Unluckly, this is not still enough. During the last meeting
devoted to Discrete Tomography, Gabor T. Herman~\cite{GTH} and
Peter Gritzmann~\cite{PG} stress the fact that various approaches
have been suggested for solving the general problem of
reconstructing binary matrices that are given by their discrete
X-rays in a small number of directions, but more work has to be
done to handle the ill-posedness of the problem. In fact, the
relevant measure for the quality of a solution of the problem
would be its deviation from the original matrix. Hence in order to
establish this deviation we would have to know the real binary
matrix. However, the goal is to find this unknown original binary
matrix so we can only consider measures for the quality of a
solution based on the given input discrete X-rays. We have a good
solution in this sense if its discrete X-rays in the prescribed
directions are close to those of the original matrix.
Unfortunately, if the input data do not uniquely determine the
matrix even a solution having the given discrete X-rays may be
very different from the unknown original matrix. It is shown in
\cite{AGT} that extremely small changes in the data may lead to
entirely different solutions. Consequently, the problem is
ill-posed, and in a strict mathematical setting we are not able to
solve this problem and get the correct solution.

In most practical application we have some a priori information
about the images that have to be reconstructed. So, we can tackle
the algorithmic challenges induced by the ill-posedness by
limiting the class of possible solutions using appropriate prior
information. The reconstruction algorithms can take advantage of
this further information to reconstruct the binary images.

A first approach is given in \cite{MVH}, where it is posed the
hypothesis that the binary matrix is a typical member of a class
of binary matrices having a certain Gibbs distribution. Then, by
using this information we can limit the class of possible
solutions to only those which are close to the given unknown
binary matrix. A modified Metropolis algorithm based on the known
Gibbs prior provides a good tool to move the reconstruction
process toward the correct solution when the discrete X-rays by
themselves are not sufficient to find such solution.

A second approach modelled a priori information in terms of a
subclass of binary images to which the solution must belong.
Several papers study the problem on classes of binary matrices
having convexity or connectivity properties. By using these
geometric properties we reduce the class of possible solutions.
For instance, there is a uniqueness result~\cite{GG} for the
subclass of {\em convex binary matrices} (i.e., finite subsets $F$
of $\ZZ^n$ such that $F=\ZZ^n \cap conv(F)$). It is proved that a
convex binary matrix is uniquely determined by its discrete X-rays in
certain prescribed sets of four directions or in any seven
non-parallel coplanar directions. Moreover, there are efficient
algorithms for reconstructing binary matrices of these subclasses
defined by convexity or connectivity properties. For example,
there are polynomial time algorithms to reconstruct {\em hv-convex
polyominoes}~\cite{BDNP,DN,BBDN} (i.e., two-dimensional binary
matrices which are 4-connected and convex in the horizontal and
vertical directions) and convex binary matrices~\cite{BD,BDD} from
their discrete X-rays. At the moment, several researchers are
studying the following stability question: given a binary matrix
having some connectivity and convexity properties and its discrete
X-rays along three or four directions, is it possible that small
changes in the data lead to ``dramatic'' change in the binary
matrix ?

In this paper, we take the second approach into consideration, and
we propose some new subclasses consisting of binary matrices with
periodicity properties. The periodicity is a natural constraint
and it has not yet been studied in Discrete Tomography. We provide
a polynomial-time algorithm for reconstructing $(p,1)$-periodical
binary matrices from their discrete X-rays in the horizontal and
vertical directions (i.e., row and column sums). The basic idea of
the algorithm is to determine a polynomial transformation of our
reconstruction problem to 2-Satisfiability problem which can be
solved in linear time~\cite{APT}. A similar approach has been
described in ~\cite{BDNP,CD}.

We wish to point out that this paper is only an initial approach
to the problem of reconstructing binary matrices having
periodicity properties from a small number of discrete X-rays.
There are many open problems on these classes of binary matrices
of interest to researchers in Discrete Tomography and related
fields: the problem of uniqueness, the problem of reconstruction
from three or more X-rays, the problem of reconstructing binary
matrices having convexity and periodicity properties, and so on.

\section{Definitions and preliminaries}
\label{s1} \setcounter{equation}{0}

\paragraph{Notations.}

Let $A^{m\times n}$ be a binary matrix, $r_i=\sum_{j=1}^n a_{i,j}$
and $c_j=\sum_{i=1}^ma_{i,j}$, for each $1\leq i \leq m$ and
$1\leq j \leq n$. We define $R=(r_1,\dots,r_m)$ and
$C=(c_1,\dots,c_n)$ as the vectors of {\em row} and {\em column
sums} of $A$, respectively. The enumeration of the rows and columns
of $A$ starts with row $1$ and column $1$ which intersect in the
upper left position of $A$. A realization of $(R,C)$ is a matrix
$B$ whose row and column sums are $R$ and $C$.

\noindent A binary matrix $A^{m\times n}$ is said to be {\em
$(p,q)$-periodical} if $a_{i,j}=1$ implies that

$a_{i+p,j+q}=1$ if $1\leq i+p \leq  m$ and $1\leq j+q \leq n$,

 $a_{i-p,j-q}=1$ if $1\leq i-p \leq  m$ and $1\leq j-q \leq n$.

\noindent
Such a matrix is said to have {\it period $(p,q)$}.

\noindent For any given couple $(x,y)$ such that $a_{x,y}=1$ we
define the set $P$ of propagation of the value in position $(x,y)$
in direction $(p,q),$ $P=\{(x+kp,y+kq)| 1\leq x+kp\leq m,1\leq
y+kq\leq n, k \in \ZZ \}$. Such set is called a {\em line}. Each
line has a {\em starting point}, which is its leftmost point, and
an {\em ending point}, which is its rightmost point. We say that a
line starts on column $j$ and ends on column $j'$ when its
starting and ending points are on column $j$ and $j'$,
respectively.

\noindent The notion of  {\em box} is a crucial part for our work. Let
$A$ be a $(p,q)$-periodical matrix. From the periodicity it
follows that if there exists an index $1\leq i \leq m$ such that
\begin{description}
\item $r_{i}=r_{i+p}+k$, then the positions on row $i$, from column $n-q+1$
to column $n$, contain at least $k$ elements equal to $1$. Such
positions form a box at the end of row $i$ and will be addressed
to as {\em right box} (rt);
\item $r_{i}+k=r_{i+p}$, then on row $i+p$, from column $1$
to column $q$ we have $k$ elements equal to $1$. Such positions
form a box at the beginning of the row $i$ and will be addressed
to as {\em left box} (lt);
\end{description}

\begin{figure}[htb]
  \centerline{\epsfxsize=5.5cm \epsfbox{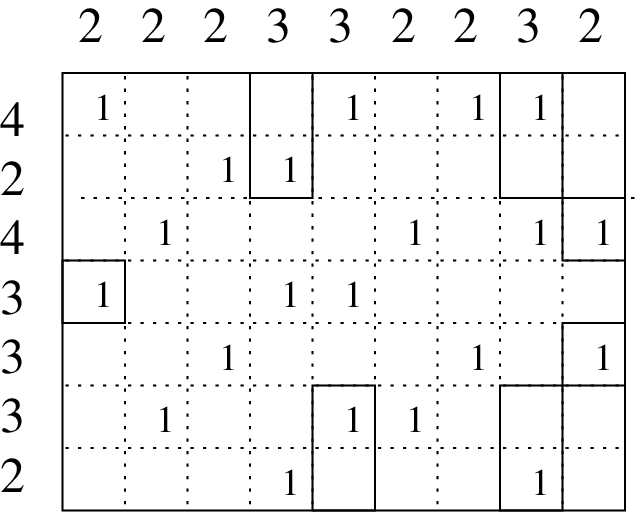}}
  \caption{\baselineskip=0ex \small
   A $(2,1)$-periodical binary matrix with upper (up), lower (lw), right
   (rt) and left (lt) boxes.}
\label{fig2}
\end{figure}

\noindent We define the upper and lower boxes (up and lw
respectively) on columns in the same way (see Fig.~\ref{fig2}), if
there exists an index $1\leq j \leq n$ such that
\begin{description}
\item $c_{j}=c_{j+q}+k$ then the positions on column $j$, from row $m-p+1$
to row $m$, contain at least $k$ elements equal to $1$. Such
positions form a box at the end of column $j$ and will be
addressed to as lower box (lw);
\item $c_{j}+k=c_{j+q}$ then the positions on column $j+q$, from row $1$ to
row $p$, contain at least $k$ elements equal to $1$. Such
positions form a box at the beginning of column $j$ and will be
addressed to as upper box (up);
\end{description}

\paragraph{Definitions of polyominoes.}

\begin{figure}[htb]
  \centerline{\epsfxsize=12cm \epsfbox{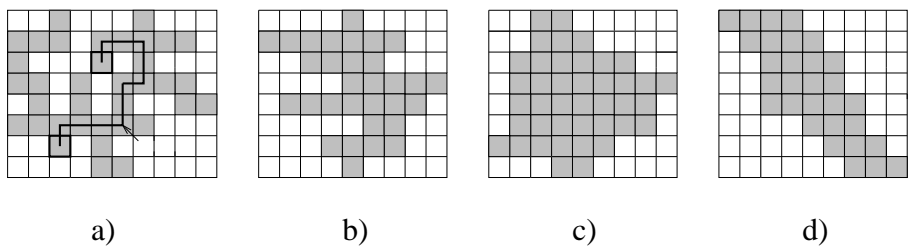} }
  \caption{\baselineskip=0ex \small
   a) A polyomino. b) A $h$-convex polyomino. c) A $hv$-convex polyomino.
   d) A parallelogram polyomino.}
\label{fig3}
\end{figure}

A {\em polyomino} $P$ is a finite union of elementary cells of the
lattice $\mathbb{Z}\times \mathbb{Z}$ whose interior is connected.
This means that, for any pair of cells of $P$ there exists a
lattice path in $P$ connecting them (see Fig.~\ref{fig3}(a)). A
lattice path is a path made up of horizontal and vertical unitary
steps. These sets are well-known combinatorial objects \cite{Gol}
and are called {\em digital 4-connected sets} in discrete geometry
and computer vision. We point out that a polyomino can be easily
represented by a binary matrix.

\noindent A polyomino is said to be {\em  $v$-convex} [{\em
$h$-convex}], when its intersection with any vertical [horizontal]
line is convex. A polyomino is {\em $hv$-convex} or {\em simply
convex} when it is both horizontal and vertical convex. A {\em
parallelogram polyomino} is a polyomino whose boundary consists of
two non intersecting paths (except at their origin and extremity)
having only north or west steps. Fig.~\ref{fig3} shows polyominoes
having the above-mentioned geometric properties.

\section{Periodicity $(1,1)$}

Let $A$ be a (1,1)-periodical matrix.

\noindent
By definition of boxes for $p=1$ and $q=1$ the boxes are reduced
to only one cell and the integer $k$ of the definition takes only
the values $0$ or $1.$ If there exists an index $1 \leq i \leq m$
or $1 \leq j \leq n$  such that
\begin{description}
\item $r_{i}=r_{i+1}+1$ then $a_{i,n}=1$.
\item $r_{i}+1=r_{i+1}$ then $a_{i+1,1}=1$.
\item $c_{j}=c_{j+1}+1$ then $a_{m,j}=1$.
\item $c_{j}+1=c_{j+1}$ then $a_{1,j+1}=1$.
\end{description}

\noindent A preprocessing part uses the previous box properties to
extract the fixed part (called $F$) of the reconstruction matrix.
The following algorithm performed on a given pair of vectors
$(R,C)$ gives, if a solution exists, the fixed part (namely the
matrix $F$) and a pair of  vectors
$(R',C')$ of the mobile part.\\

\noindent
{\bf Propagation}$(x,y,F,\ell)$
\begin{tabbing}
$P=\{(x+k,y+k)| 1\leq x+k\leq m,1\leq y+k\leq n, k \in \ZZ \};$\\
{\bf  For all} $(i,j) \in P$ {\bf do} $F_{i,j}=1,
r_i^{(\ell+1)}=r_i^{(\ell)}-1, c_j^{(\ell+1)}=c_j^{(\ell)}-1;$
\vspace{0.3cm}
\end{tabbing}

The main program finds the fixed 1's by considering the
differences between the values of the pair of vectors. For
each fixed 1, the procedure Propagation fills by periodicity
$(1,1)$ the matrix $F$ and decreases the row and column sum of the
current matrix.

\noindent At the end of the preprocessing part both row and
column vectors are homogeneous ($r'_i=\rho$ for all $1 \leq i \leq
m$ and $c'_j=\gamma$ for all $1 \leq j \leq n$).

\noindent Now, either  $(R',C')$ are zero vectors and then the
solution is unique and equal to $F$, or we perform a
reconstruction from homogeneous X-rays with periodicity $(1,1)$.

Since the vectors are homogeneous we can extend the periodicity on
a torus. Indeed, suppose that $r'_i=r'_{i+1}$ with $r'_i=
\sum_{j=1}^n a'_{i,j}$ and $r'_{i+1}= \sum_{j=1}^n a'_{i+1,j}$. By
periodicity $a'_{i,j}= a'_{i+1,j+1}$, for $j=1, \ldots, n-1$,
implies that $ a'_{i,n}= a'_{i+1,1}$, for $i=1, \ldots, m-1$. In
other terms the values of the matrix $A'$ are mapped on a
cylinder. The same argument in column proves that the values of
the matrix $A'$ are mapped on a torus. That is if $a'_{i,j}=1$
then $a'_{i+1 \mod m ,j+1 \mod  n}=1$ and $a'_{i-1 \mod m ,j-1
\mod n}=1$.

\noindent So, a solution is formed by loops, namely a beginning in
$(i,j)$ with $a'_{i,j}=1$ and a propagation by periodicity $(1,1)$
until the position $(i,j).$ All the loops have the same length. As
the vectors are homogeneous, we can compute the number of loops.
Using this strong condition and the algorithm of Ryser~\cite{Ry}
in the first row in order to place the loops, we can reconstruct
easily a solution in $O(m n)$ time.

Another remark is the arithmetical nature of the stability of the
solution. We can prove that if $m$ and $n$ are relatively prime
then there is only one solution. Indeed in this case, to perform a
reconstruction of a binary matrix with homogenous vector, we have
only one loop felling the whole matrix and then either the matrix
$A'$ is full of 1's or full of 0's and nothing between because of
the toric conditions.
\begin{proposition}
Let $R\in \IN^m$ and $C\in \IN^n$. If $\gcd(n,m)=1$, then
there is at most a $(1,1)$-periodical matrix having row and column sums equal
to $(R,C)$.
\end{proposition}
For example, if we perform a reconstruction with a matrix with $m$
rows and $m+1$ columns and periodicity $(1,1)$ then the solution
is unique if it exists.

\begin{example}
For $R=(2,2,1,2), C=(2,1,2,2)$, the algorithm gives the matrix

$F=\begin{array}{llll}
0 & 0 & 1&0 \\
0 & 0 & 0&1 \\
0 & 0 & 0&0 \\
1 & 0 & 0&0 \\
\end{array}$\\
and the vectors $R=(1,1,1,1), C=(1,1,1,1)$.\\
We can reconstruct two solutions for $R=(2,2,1,2), C=(2,1,2,2)$:

$A'_1=\begin{array}{llll}
1 & 0 & 1&0 \\
0 & 1 & 0&1 \\
0 & 0 & 1&0 \\
1 & 0 & 0&1  \\
\end{array}$\\
and

$A'_2=\begin{array}{llll}
0 & 0 & 1&1 \\
1 & 0 & 0&1 \\
0 & 1 & 0&0 \\
1 & 0 & 1&0 \\
\end{array}.$\\

For $R=(1,2,2,2), C=(2,1,2,1,1)$, the algorithm gives a fixed part
and $R=(0,0,0,0), C=(0,0,0,0,0)$, then the solution is unique.

$F=\begin{array}{lllll}
0 & 0 & 1&0 &0 \\
1 & 0 & 0&1 &0\\
0 & 1 & 0&0 &1\\
1 & 0 & 1&0 &0\\
\end{array}$.
\end{example}

\section{Periodicity $(p,1)$ with $ 1<p<m$}

Let $A$ be a matrix with periodicity $(p,1).$

\noindent
The preprocessing part uses only the row sums in order to find the
fixed part of the reconstruction. In fact, by definition of boxes for
$q=1$ the horizontal boxes are reduced to only one cell and the
integer $k$ of the definition takes only the values $0$ or $1.$
If there exists an index $1 \leq i \leq m$
such that
\begin{description}
\item $r_{i}=r_{i+p}+1$ then $a_{i,n}=1$.
\item $r_{i}+1=r_{i+1}$ then $a_{i+p,1}=1$.
\end{description}

A preprocessing part uses the previous box properties to extract
the fixed part (called $F$) of the reconstruction matrix. The
following algorithm performed on a given pair of vectors $(R,C)$
gives, if a solution exists, the fixed part (namely the matrix
$F$) and a pair of $(R',C')$ of the mobile part.\\

\begin{tabbing}
{\bf Algorithm 2}\\ \\
{\bf Input:}  A pair of integral vectors $(R,C)$; \\
\= {\bf Output:}\= If $PB=0$, then it gives Matrix F and couple of integral vector
$(R',C')$;\\
 \> \>{\bf or}\ If $PB=1$ Failure in the reconstruction;\\
{\bf  Initialisation:} $\ell:=0, PB:=0, R^{(\ell)}:=R,C^{(\ell)}:=C, F=0^{m\times n}$;\\

\= {\bf while} \= $R^{(\ell)}$ is positive non homogeneous vector  and $PB=0$ {\bf do}\\
\> \> $R^{(\ell+1)}:=R^{(\ell)},C^{(\ell+1)}:=C^{(\ell)};$\\
\> \> {\bf determine} first index s.t. $r_i^{(\ell)} \neq r_{i+p}^{(\ell)}$; \\
\> \> {\bf if} $r_{i}^{(\ell)}=r_{i+p}^{(\ell)}+1$ {\bf then}
$x:=i,y:=n,$
{\bf Propagation}$(x,y,F,\ell)$;\\
\> \> {\bf else} \=  \={\bf if}
$r_{i}^{(\ell)}+1=r_{i+p}^{(\ell)}$ {\bf then} $x:=i+p,y:=1,$
{\bf Propagation}$(x,y,F,\ell)$\\
\> \> \> \>{\bf else}  $PB:=1;$\\
\> \> $\ell:=\ell+1;$\\
{\bf end while};\\
\vspace{0.3cm}
\end{tabbing}

\noindent
{\bf Propagation}$(x,y,F,\ell)$
\begin{tabbing}
$P=\{(x+kp,y+k)| 1\leq x+kp\leq m,1\leq y+k\leq n, k \in \ZZ \};$\\
{\bf  For all} $(i,j) \in P$ {\bf do} $F_{i,j}=1,
r_i^{(\ell+1)}=r_i^{(\ell)}-1, c_j^{(\ell+1)}=c_j^{(\ell)}-1;$
\vspace{0.3cm}
\end{tabbing}

The main program finds the fixed 1's by considering the
differences between the values of the row sums. For each fixed 1,
the procedure Propagation fills by periodicity $(p,1)$ the matrix
$F$ and decreases the row and column sum of the current matrix.

\noindent At the end of the preprocessing part the row  vector sum
$R'$ of $A'$ has the same value on indices in arithmetical
progression of rank
$p$: $r'_i=r'_{i+p}=r'_{i+2p} \cdots
=r'_{i+(\ell-1) p}$ where $\ell= L $ or $L+1$. This set of element
of the row sums $R'$ of $A'$ is called {\em line} of $R'$. The
minimum length of each line of $R'$ is $L=\lfloor \frac{m}{p}
\rfloor$. The number of lines of length $L+1$ and $L$ of $R'$ is
$n_{L+1}=m \mod p$ and $n_L=p-n_{L+1}$, respectively.

\begin{example}
If  $(p,1)=(2,1)$, and $R=(2,3,2,4,3,4,2), C=(3,4,3,3,4,2,1)$, the algorithm gives the matrix

$F=\begin{array}{lllllll}
0 & 0 & 0 & 0 & 1 & 0 &0 \\
0 & 0 & 0 & 0 & 0 & 0 &0 \\
0 & 0 & 0 & 0 & 0 & 1 & 0 \\
1 & 0 & 0 & 0 & 0 & 0 & 0 \\
1 & 0 & 0 & 0 & 0 & 0 & 1 \\
0 & 1 & 0 & 0 & 0 & 0 & 0 \\
0 & 1 & 0 & 0 & 0 & 0 & 0 \\
\end{array}$\\
and the new vectors are $R'=(1,3,1,3,1,3,1), C'=(1,2,3,3,3,1,0)$.
Since $m=7$ and $p=2$, we have that $L=3, n_L=1, n_{L+1}=1$.
So, $R'=(1,3,1,3,1,3,1)$ contains a line of length $L=3$ and
a line of length $L+1=4$.
Since the lines of length $L=3$ and $L+1$ are $(*,3,*,3,*,3,*)$ and
$(1,*,1,*,1,*,1)$, respectively, we have that matrix $A'$ contains
a line of length
three lines of length $L=3$ starting from the second row and
$L+1=4$ starting from the first row.
\end{example}

Now, we prove now the values of $A'$ are mapped on a cylinder. We
have  $r'_i=r'_{i+p}$ with $r'_i= \sum_{j=1}^n a'_{i,j}$ and
$r'_{i+p}= \sum_{j=1}^n a'_{i+p,j}$. From the periodicity
$a'_{i,j}= a'_{i+p,j+1}$, for $j=1, \cdots, n-1$, it follows that
$ a'_{i+p,1}= a'_{i,n}$, for $i=1, \cdots m-p$. In other terms the
values of the matrix $A'$ are mapped on a cylinder.

Thus a 1 on the first  $p$ rows (at position $(x,y), 1\leq x \leq
p, 1 \leq y \leq n )$ can be extended by periodicity on the matrix
$A'$ by 1's in positions $(x + kp ,y +k \mod n)$ with $k= 0,
\cdots, \ell-1$ where $\ell= L $ or $L+1$. The matrix $A'$ is in
particular composed on a cylinder of lines in direction $(p,1)$ of
length $L$ or $L+1$. In addition to that the number of lines of
length $L+1$ is exactly $n_{L+1}=m  \mod  p$ and the number of
lines of length $L$ is $n_L=p- n_{L+1}$.

\subsection{A reduction to the problem of reconstructing
a special class of h-convex binary matrices lying on a cylinder}
\label{ss}

Let $A'$ be a solution for a given $(R',C')$ with
$r'_i=r'_{i+p}=r'_{i+2p} \cdots =r'_{i+(\ell-1)}$ where $\ell= L $
or $L+1$. We now perform a reduction of reconstruction of
$(p,1)$-periodical matrix $A'$ on a cylinder to a reconstruction
of a special $h$-convex matrix $A''$ on a cylinder.

By the previous construction, matrix $A'$ is formed by lines $(x +
kp ,j +k \mod n)$ with $k= 0, \ldots, \ell-1$ where $\ell= L$ or
$L+1$. The starting points of the lines is the set of position
$S=\{(x,y)| 1 \leq x \leq p, 1 \leq y \leq n, a'_{x,y}=1\}$.

\noindent $S$ is ordered by: $(x,y) \leq (x',y')$
if and only if $x \geq x'$ and $y \leq y'$ (i.e., we proceed from
bottom to up and from left to right).

\noindent Let $S'$ be the set $S$ with an extra index of the rank
in the previous order. Each element of $S'$ is a triple $(x,y,o)$,
where $(x,y)$ is an element of $S$ and $o$ the rank in the order.
Now, we can describe the reduction.\\

\noindent {\bf Reduction.} Let $(x,y,o)\in S'$. The point $(x,y)$
is the starting point of a line of matrix $A'$ having length $L$
or $L+1$. This line gives a horizontal bar of 1's begins in
position $(o ,y +r \mod n)$ with $r= 1, \ldots, \ell$ where $\ell=
L$ if $x > n_{L+1}$ or $\ell=L+1$ if $x \leq n_{L+1}$. The set of
these horizontal bars gives the $h$-convex matrix $A''$.

\

\noindent Notice that,  this  transformation makes the column sum
of $A''$ equal to the column sum of $A'$. In Section~\ref{sss}, we
will show the inverse reduction that provides $A'$ from $A''$.

\begin{example}
Let us take the following matrix $A'$ with periodicity $(2,1)$ and
$R'=(3,2,3,2,3)$ and $C'=(3,3,2,1,2,2)$ into consideration.

$A'=\begin{array}{llllll}
1 & 0 & 1 & 0 & 0 &1\\
1 & 0 & 0 & 0 & 1 &0\\
1 & 1 & 0 & 1 & 0 &0 \\
0 & 1 & 0 & 0 & 0 &1\\
0 & 1 & 1 & 0 & 1 &0\\
\end{array}$\\
Since $m=5$ and $p=2$, we have that $L=2$, $n_{L+1}=1$. The matrix
$A'$ is composed of three lines of  length $L+1=3$ and two lines
of length $L=2$. The starting points in the first two rows (the
two first indices are the position in the matrix $A'$ and the last
index is the rank in the order) are:
$$S'=\{(2,1,1),(1,1,2),(1,3,3),(2,5,4),(1,6,5)\}.$$
The transformation gives the following  $h$-convex matrix $A''$
mapped on a cylinder with three bars of length three and two bars
of length two.

$A''=\begin{array}{llllll}
1 & 1  & 0 & 0 & 0 & 0\\
1 & 1 & 1 & 0 & 0 & 0\\
0 & 0 & 1 & 1 & 1 & 0 \\
0 & 0 & 0 & 0 & 1 & 1\\
1 & 1 & 0 & 0 & 0 & 1\\
\end{array}$\\
The column sums of $A''$ are equal to $C'=(3,3,2,1,2,2)$.
\end{example}
We point out that the order on the starting points adds the
following constraints:\\

\begin{figure}[htb]
  \centerline{\epsfxsize=4cm \epsfbox{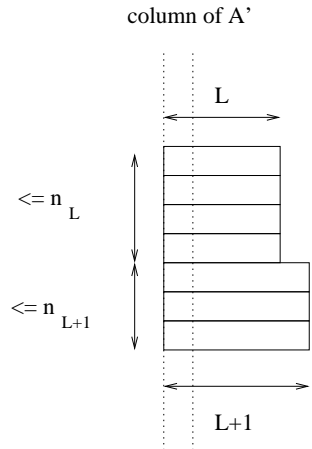} }
  \caption{\baselineskip=0ex \small
   \label{fig4}A column of a binary matrix satisfying
condition 1.}
\end{figure}

\noindent {\bf Condition 1.} On each column of $A''$ can start at
most $n_L$ bars of length $L$ and at most $n_{L+1}$ bars of length
$L+1$. Moreover, by proceeding from up to down on the column at
first we find the bars of the length $L$ and then the bars of the
length $L+1$ (see Fig.~\ref{fig4}).

\

\noindent We denote the class of $h$-convex binary matrices lying
on a cylinder and satisfying condition~1 by
$\mathcal{HC}(n_{L},n_{L+1})$.

\noindent By this property, the matrix $A''$ consists of four
disjoint zones $B, C, E$ and $P$ whose boundaries are three paths
having only north or west steps, and $a''(i,j)=1$ for $(i,j)\in
C\cup P$ and $a''(i,j)=0$ for $(i,j)\in B\cup E$ (see
Fig.~\ref{fig5} and the matrix $A''$ of the previous example).
Notice that, the matrices of the class
$\mathcal{HC}(n_{L},n_{L+1})$ are set of parallelogram polyominoes
lying on a cylinder.

\

From the reduction it follows that the problem of reconstructing a
$(p,1)$-periodical binary matrix $A'$ having row and column sums
$(R',C')$ (output of Algorithm 1) is equivalent to the problem of
reconstructing a binary matrix  $A''$ of
$\mathcal{HP}(n_{L},n_{L+1})$ having column sums $C'$, $m$ rows of
length $L$ and $L+1$. We denote this reconstruction problem on the
cylinder by ${\bf RHC}$ problem.

In the following subsection, we determine a polynomial
transformation of ${\bf RHC}$ problem to 2-Satisfiability problem
(2-SAT).

\subsection{A reduction to the 2-SAT problem}

Given an instance $I$ of  ${\bf RHC}$ problem, we want to build a
2-SAT formula $\Omega$ (a formula in conjunctive normal form,
where each clause has at most two literals) whose satisfiability
is linked to the existence of a solution for $I$ in such a way: if
$\Omega$ is satisfiable, then we are able to reconstruct a
solution for $I$ in P-time and, vice versa, each solution of $I$
gives  an evaluation of the variables satisfying $\Omega$ in
P-time. We will do not show the proofs of the lemmas of this
section for brevity's sake. Let  $I$ be an instance of ${\bf RHC}$
problem; that is:
\begin{itemize}
\item[-] two integers $n_L, n_{L+1}$;
\item[-] a couple $(L,C)$, where $L$ and $L+1$ are the only possible
values of the row sums of a binary matrix $A''$ of
$\mathcal{HC}(n_L, n_{L+1})$, solution of $I$, and
$C=(c_1,\dots,c_n)$ is its column sums.
\item[-] an integer $m$ which denotes the number of rows of $A''$.
\end{itemize}

The formula $\Omega$ that we want to construct is the conjunction
of three $2$-SAT formulas: $\Omega_1$ which encodes the
geometrical constraints of $A''$, $\Omega_2$ which gives the
consistency of $A''$ with the couple $(L,C)$ and, finally,
$\Omega_3$ which imposes the constraints of condition 1 on each
column of $A''$. The variables of the formula $\Omega$ belong to
the union of the four disjoint sets of variables:
$$\mathcal{B}=\{b(i,j):1\leq i \leq m, 1\leq j \leq n\} \mbox{ , } \:\:\:\:\:
\:\:\:\:\: \mathcal{C}=\{c(i,j):1\leq i \leq m, 1\leq j \leq n\}
\mbox{ , }
$$
$$\mathcal{P}=\{p(i,j):1\leq i \leq m, 1\leq j \leq n\} \mbox{ , } \:\:\:\:\:
\:\:\:\:\: \mathcal{E}=\{e(i,j):1\leq i \leq m, 1\leq j \leq
n\}.$$
We use the variables of the set $\mathcal{X}_i$, with
$\mathcal{X}_i \in \{\mathcal{B}, \mathcal{C}, \mathcal{P},
\mathcal{E}\}$, to represent the four disjoint zones $B$, $C$,
$P$, $E$ inside $A$.

\subsubsection{Coding in $\Omega_1$ the geometrical constraints of $A''$.}

Formula $\Omega_1$ is the conjunction of the following sets of
clauses:
$$
\begin{array}{llll}
Corners & = & \bigwedge_{i,j}  \left\{          
\begin{array}{l}
(x(i,j)\Rightarrow x(i-1,j)) \:\wedge \: (x(i,j)\Rightarrow
x(i,j+1)) \:\: \mbox{for $x\in \mathcal{C}\bigcup \mathcal{E}$} \\ \\
(x(i,j)\Rightarrow x(i+1,j))\:\wedge\: (x(i,j) \Rightarrow
x(i,j-1)) \:\: \mbox{for $x\in \mathcal{B}\bigcup \mathcal{P}$}
\end{array}
\right\} \\
\\
Disj & = & \bigwedge_{i,j}  \left\{  
\begin{array}{l}                    
(b(i,j)\Rightarrow \overline{c}(i,j)) \:\wedge \: (p(i,j)
\Rightarrow b(i,j)) \:\wedge \: (e(i,j) \Rightarrow c(i,j))
\end{array}
\right\}\\
\\
Compl &= & \bigwedge_{i,j}  \{  
\overline{b}(i,j)\Rightarrow c(i,j) \}\\ \\
Anch &=& \{\overline{e}(1,L) \wedge \overline{p}(m,L+1)\wedge p(r,1)\}
\end{array}
$$
with $ 1< r \leq m$.
\begin{definition}
Let $V_1$ be an evaluation of the variables in $\mathcal{B}$,
$\mathcal{C}$, $\mathcal{P}$, $\mathcal{E}$ which satisfies
$\Omega_1$. We define the binary matrix $A''$ of size $m\times n$
as follows:
$$
\begin{array}{ll}
(c(i,j)=1 \;\; \wedge \;\; e(i,j)=0) \; \Rightarrow \; a''(i,j)=1
\mbox{ ,} \qquad p(i,j)=1 \; \Rightarrow \; a''(i,j)=1 \mbox{ ,} \\
(b(i,j)=1 \;\; \wedge \;\; p(i,j)=0) \; \Rightarrow \; a''(i,j)=0
\mbox{ ,} \qquad e(i,j)=0 \; \Rightarrow \; a''(i,j)=0.
\end{array}
$$
\label{def1}
\end{definition}
It is immediate to check that $A''$ is well defined.

\begin{figure}[htb]
  \centerline{\epsfxsize=6.5cm \epsfbox{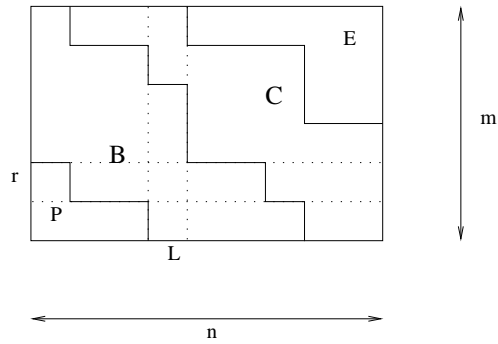} }
  \caption{\baselineskip=0ex \small
   \label{fig5}The matrix $A''$ and the four zones $B, C, E$
and $P$ of $A$ defined by the formulas $Corners, Disj, Compl$ and
$Anch$.}
\end{figure}

\noindent The matrix $A''$ contains the four zones $B, C, E$ and
$P$ of $A$ such that: $(i,j)\in X$, with $X\in \{B,C,E,P\}$ if and
only if $x(i,j)=1$ (see Fig.~\ref{fig5}). From $Corners$, $Disj$,
$Compl$ and $Anch$ we deduce the following properties these four
zones:
\begin{lemma}\label{lem1}
\begin{itemize}
\item[i)] $\{B,C\}$ is a partition of $A''$, $P\subseteq B$ and $E\subseteq
C$; \item[ii)] the boundary of zones $B, C, E$ and $P$ is made up
of three paths having only north or west steps;
\item[iii)]
there does not exist a column of $A''$ containing both points of
$P$ and points of $E$.
\end{itemize}
\end{lemma}

\subsubsection{Coding in $\Omega_2$ the bound of the
row and column sums of $A''$.}

The formula $\Omega_2$ is the conjunction of the following sets of
clauses:
$$
\begin{array}{l}
LBC = \bigwedge_{i,j} \left\{
\begin{array}{ll}
\mbox{ if } j>L, & e(i,j)\Rightarrow \overline{b}(i+c_j,j)\\
\\
\mbox{ if } j\leq L, & b(i,j) \Rightarrow p(i+m-c_j,j)
\end{array}
\right\} \\
\\
UBC = \bigwedge_{i,j} \left\{
\begin{array}{ll}
\mbox{ if } j>L, & \overline{e}(i,j) \Rightarrow b(i+c_j,j) \\
\\
\mbox{ if } j\leq L, & \overline{b}(i,j) \Rightarrow
\overline{p}(i+m-c_j,j)
\end{array}
\right\}\\
\\
UBR = \bigwedge_{i,j} \left\{
\begin{array}{ll}
\mbox{ if } i<r, & \overline{b}(i,j) \Rightarrow e(i,j+L+1)\\
\\
\mbox{ if } i\geq r, & p(i,j) \Rightarrow \overline{c}(i,j+n-L-1)
\end{array}
\right\}\\
\\
LBR = \bigwedge_{i,j} \left\{
\begin{array}{ll}
\mbox{ if } i<r, & b(i,j) \Rightarrow \overline{e}(i,j+L)\\
\\
\mbox{ if } i\geq r, & \overline{p}(i,j) \Rightarrow c(i,j+n-L)
\end{array}
\right\}.
\end{array}
$$
The formulas $LBC$ and $UBC$ give a lower and an upper bound
for the column sums of $A$.
The formula $LBR$ express that the row sums are greater than $L$.
Finally, the formula $UBR$ express that the row sums are smaller than $L+1$.
More precisely,
\begin{lemma}\label{lem2}
Let $A''$ be the binary matrix defined by means of the valuation
$V_2$ which satisfies $\Omega_1 \wedge \Omega_2$ as in
Definition~\ref{def1}. We have that:
\begin{description}
\item{$i)$} the column sums of $A''$ are equal to $C=(c_1,\dots,c_n)$;
\item{$ii)$} each row sum of $A''$ has value $L$ or $L+1$ (see Figure
\ref{fig5}).
\end{description}
\end{lemma}

\begin{figure}[htb]
  \centerline{\epsfxsize=10.5cm \epsfbox{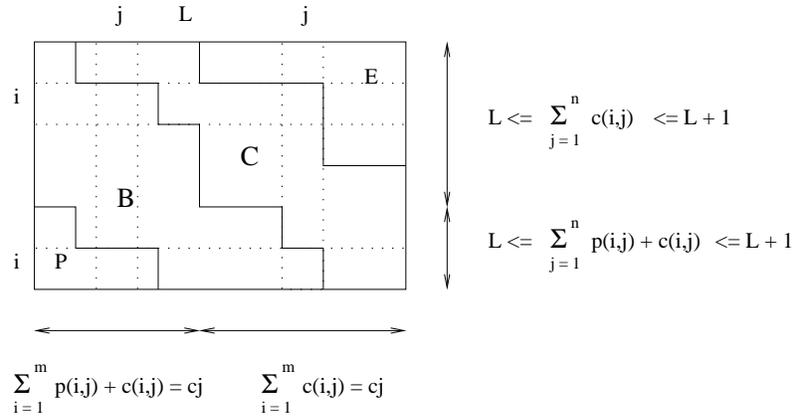} }
  \caption{\baselineskip=0ex \small
   \label{fig6}The matrix $A''$ with the bounds on row and
column sums.}
\end{figure}

\subsubsection{Coding in $\Omega_3$ the maximum number of
bars of length $L$ and $L+1$ starting on each column of $A''$.}

The formula $\Omega_3$ is the conjunction of the following sets of
clauses:
$$
\begin{array}{ll}
BB_{L} &= \bigwedge_{i,j}\left\{
\begin{array}{ll}
\mbox{ if } i\leq m-c_n, & \:\:  e(i,j) \Rightarrow
\overline{b}(i-n_{L+1},j-L-1)\\
\\
\mbox{ if } m-c_n<i<r, & \:\: r-(m-c_n)-1>n_{L+1} \Rightarrow
\overline{b}(m-c_n+1,n-L)\\
\\
\mbox{ if } (i\geq r \wedge j\leq L), & \:\: \overline{p}(i,j)
\Rightarrow \overline{b}(i-n_{L+1},j+n-L-1)
\end{array} \right\}\\
\\
BB_{L+1} & = \bigwedge_{i,j}\left\{
\begin{array}{ll}
\mbox{ if } i\leq n-c_n, & \:\: c(i,j) \Rightarrow e(i-n_L,j+L) \\
\\
\mbox{ if } n-c_n <i<r, & \:\: r-(m-c_n)-1 > z_1 \Rightarrow
\overline{c}(m-c_n+n_L+1,n-L)\\
\\
\mbox{ if } i\geq r, & \:\: c(i,j) \Rightarrow
\overline{p}(i-n_L,j+L-n)
\end{array}
\right\}
\end{array}
$$

\begin{lemma}\label{lem3}
Let $A''$ be the binary matrix defined by means of the valuation
$V_3$ which satisfies $\Omega_1 \wedge \Omega_2 \wedge \Omega_3$
as in Definition~\ref{def1}. We have that:
\begin{description}
\item{$i)$} on each column of $A''$ can start at most $n_{L}$
bars of length $L$,
\item{$ii)$} on each column of $A''$ can start at most $n_{L+1}$
bars of length $L+1$.
\end{description}
\end{lemma}

\noindent By Lemmas~\ref{lem1} and \ref{lem3}, matrix $A''\in
\mathcal{HP}(n_{L},n_{L+1})$. By Lemma~\ref{lem2}, matrix $A''$
satisfies the tomographic constraints. Therefore:
\begin{theorem}
\label{maintheor}
  $\Omega_1 \wedge \Omega_2 \wedge \Omega_3$ is satisfiable if and only if
  there is a binary matrix $A''$ of $\mathcal{HP}(n_{L},n_{L+1})$ having column
  sum $C'$, $m$ rows of length $L$ and $L+1$.
\end{theorem}

Since $\Omega_1 \wedge \Omega_2 \wedge \Omega_3$ is a boolean
formula in conjunctive normal form with at most two literals in
each clause, by Theorem~\ref{maintheor} we have a polynomial time
transformation of {\bf HRC} problem to 2-SAT problem which can be
solved in linear time~\cite{APT}.

\subsection{Final Step}
\label{sss}

By performing the previous reduction and an algorithm for solving
2-SAT problem~\cite{APT}, we obtain a matrix $A''$ of
$\mathcal{HP}(n_{L},n_{L+1})$ having column sums $C'$, $m$ rows of
length $L$ and $L+1$, where $n_{L+1}=m  \mod  p$ and
$n_L=p-n_{L+1}$. Now, for determining a $(p,1)$-periodical matrix
$A'$ having row and column sums equal to $(R',C')$, we have to
perform the inverse of the reduction defined in Section~\ref{ss}.
We point out that $(R',C')$ is the output of Algorithm 2. This
inverse reduction should provides $A'$ from $A''$. The following
algorithm describes this inverse reduction.

\begin{tabbing}
{\bf Algorithm 3}\\ \\
\= {\bf Input:}\=  \ the matrix $A''$ whose column sums is vector
$C'$, and vector $R'$\\
\>  \> which is homogeneous with respect to $p$ ($R'$ is the output of Algorithm 2); \\
{\bf Output:} \ the matrix $A'$ having row and column sums $(R',C')$;\\

\= {\bf Step 1:} \=  {\bf determine} the two vectors
$C^{(1)}=(c^{(1)}_1,\ldots,c^{(1)}_n)$, $C^{(2)}=(c^{(2)}_1,\ldots, c^{(2)}_n)$\\
\> \> such that $c^{(1)}_i$ and $c^{(2)}_i$ are the number of bars
starting from  the $i$-th \\
\> \>column of $A''$ having length $L$ and $L+1$, respectively, ;\\

\= {\bf Step $2$:}\= \ {\bf construct} the matrix $A'$ of size $m
\times n$ in such a way:\\

\> \> \= {\bf if }\= \  $1\leq i \leq (m-p)$ and $1\leq j \leq n$,
we set $a'(i,j)=0$;\\

\> \> \= {\bf if }\= \ $(m-p+1) \leq i \leq (m-p+n_{L+1})$ and
$1\leq
j \leq n$, {\bf then} perform Ryser's\\
\> \> reconstruction algorithm on the row and column vectors:\\
\> \> $(r_{m-p+1},\ldots , r_{m-p+n_{L+1}})$ and $C^{(1)}$; \\

\> \> \= {\bf if }\= \  $(m-p+n_{L+1}+1) \leq i \leq m$ and $1\leq
j
\leq n$, {\bf then} perform Ryser's\\
\> \> reconstruction algorithm on the row and column vectors:\\
\> \> $(r_{m-p+n_{L+1}+1},\dots , r_m)$ and $C^{(2)}$;\\

\= {\bf Step $3$:}\= \ {\bf for all} $m-p+1 \leq i\leq m$ and
$1\leq j \leq n$, \\
\> \>\ {\bf if} $a'(i,j)=1$ {\bf then} $a'(i-pk,(j+k)\mod n)=1$
with $1 \leq i-pk \leq m-p$.

\vspace{0.3cm}
\end{tabbing}

\begin{proposition}
The matrix $A'$ which is the output of the Algorithm $3$ is
$(p,1)$-periodical and it has row and column sums equal to
$(R',C')$.
\end{proposition}
We do not show the proof of this Proposition for brevity's sake.

\noindent By performing Algorithm 1, the reduction of
Section~\ref{ss}, the reduction to 2-SAT problem, a linear-time
algorithm for solving 2-SAT problem~\cite{APT}, and Algorithm 3,
we obtain a $(p,1)$-periodical binary matrix $A$ having row and
column sums equal to $(R,C)$. Since each step can be perfomed in
polynomial time, we have that:
\begin{theorem}
The problem of reconstructing $(p,1)$-periodical binary matrices
from their row and column sums can be solved in polynomial time.
\end{theorem}

\section{Conclusions}

Our main purpose has been to introduce periodicity properties in
terms relevant for Discrete Tomography. The periodicity is a
natural constraint and it has not yet been studied in this field.
The motivation of this study is in the attempt to tackle the
ill-posedness of the reconstruction problem by limiting the class
of possible solutions using appropriate prior information. This
means that, we modelled a priori information in terms of a
subclass of binary images to which the solution must belong.

By using the periodicity properties we reduce the class of
possible solutions. For instance, we proved a uniqueness result
for the class of binary matrices having period $(1,1)$. We have
shown a simple greedy algorithm for reconstructing this class of
matrices from their row and column sums. This reconstruction
problem becomes more difficult for the binary matrices having
period $(p,1)$ or $(1,q)$. We have described a polynomial-time
algorithm for solving this problem which use a reduction to
2-Satisfiability problem. We stress the fact that an interesting
property of this approach is that it can be used for
reconstructing parallegram polyominoes lying a cylinder from row
and column sums.

The future challenges concern the reconstruction of binary
matrices with a generical period $(p,q)$. We wish to point out
that this paper is only an initial approach to the problem of
reconstructing binary matrices having periodicity properties from
a small number of discrete X-rays. Lot of work should be done to
understand such environment: we only challenge the reconstruction
problem from two X-rays in some special cases, but MANY
consistency, reconstruction and uniqueness problems can be
reformulated imposing periodical constraints.


\end{document}